\documentclass[10pt]{article}
 \usepackage{amsmath, amsthm, amssymb}
\usepackage{units}
\usepackage{tabularx}
\usepackage{float}
\usepackage{mathrsfs}
\usepackage[utf8]{inputenc}
\usepackage{graphics} 
\usepackage{epsfig} 
\usepackage[section]{placeins}

\newcommand{\dd}{\mathrm{d}}
\newcommand{\bT}{\textbf{T}}
\newcommand{\e}{\operatorname{e}}

\newtheorem{thm}{Theorem}[section]

\newtheorem{lemma}[thm]{Lemma}

\newtheorem{definition}[thm]{Definition}

\newtheorem{example}[thm]{Example}

\allowdisplaybreaks
\title{Stochastic Perturbation of the Lighthill-Whitham-Richards Model via the Method of Stochastic Characteristics}
\author{
	Nora M\"uller\\
	{\small Fakultät für Mathematik, Bielefeld University,}\\ 
	{\small	Universitätsstraße 25, 33615 Bielefeld, Germany}\\
	{\small E-Mail:nmueller@math.uni-bielefeld.de}\\[.3cm]
	Wolfgang Bock\\
	{\small Technomathematics Group}\\
	{\small	University of Kaiserslautern}\\
	{\small	P.\ O.\ Box 3049, 67653 Kaiserslautern, Germany}\\
	{\small E-Mail:bock@mathemaik.uni-kl.de}\\[.3cm]}

\begin{document}
\maketitle

\begin{abstract}
In this paper we apply the method of stochastic characteristics to a Lighthill-Whitham-Richards model. The stochastic perturbation can be seen as errors in measurement of the traffic density. For concrete examples we solve the equation perturbed by a standard Brownian motion and the geometric Brownian motion without drift.
\end{abstract}




\section{Introduction}

Many traffic flow models go back to scalar conservation laws of generally non-linear type, i.e.
\begin{equation}\label{scalar-cons-law}
	u_t + f(u)_x =0,
\end{equation}
see e.g.~\cite{Bressan}. The function $u$ describes the density of vehicles on a road and thus has values on the compact set $[0,1]$. A conservation law is derived under the assumption, that the time propagation of a mass on a certain interval is only affected by flux at the boundary of the interval. Hence one often chooses 
$
f(u)=u\cdot v(u),
$
where $v$ is the Eulerian velocity of the traffic. One of the most easy choices is the Lighthill-Whitham-Richards model, which uses
$f(u)=u\cdot(1-u)$, i.e.~a velocity depending linearly on the density.
The scalar conservation law (\ref{scalar-cons-law}) now reads
\begin{equation}\label{traffic-eq}
	u_t + (1-2u)\cdot u_x =0.
\end{equation}
The flux function in the Lighthill-Whitham-Richards model is in a relatively good agreement with traffic measurements, see \cite{LW55}. The main problem is that measurements show that data points are quite accurate for low and high densities but are noisy around the maximum point. The flux function is hence rather given by
$$f(u)= u\cdot (1-u) + H(u) \circ \frac{\dd M_t}{\dd t},$$
where $H$ is a function vanishing at $u=0$ and $u=1$ and $M_t$ is a suitable nice enough stochastic process.
Plugging this in the above conservation law (\ref{scalar-cons-law}) yields
\begin{equation}\label{traffic-eq+H}
	u_t + f(u)_x =u_t +(1-2u)\cdot u_x+ H(u)_x \circ \frac{\dd M_t}{\dd t}=0.
\end{equation}
In this manuscript, we use a heuristic approach of stochastic characteristics to solve such equations for the stochastic perturbed Lighthill-Whitham-Richards model. Indeed we will work out explicit solutions for different cases, where the driving process is given by a Brownian motion or a geometrical Brownian motion. This work can be seen as the starting point to the investigation of different stochastically perturbed hyperbolic equations. 

\section{Prelimiaries}
We look at some examples for the Lighthill-Whitham-Richards model for different initial conditions. Due to the underlying model the initial condition describes the density of our traffic problem at time $t=0$ and at position $x\in[0,1]$. 

\begin{lemma}\label{LE:H=0_CE}
	Consider the following partial differential equation on $[0,1]$
	\begin{equation}\label{eq:LWR_H=0}
		\left\{\begin{aligned}
			\dd u &=-(1-2u)\cdot u_x\,\dd t,\\
			u(x,0)&= g(x),
		\end{aligned}\right.
	\end{equation}
	where $g(x)$ is a smooth function. Let $(\xi_t,\eta_t)$ be the solutions to the so called characteristic equations given by
	\[
	\left\{\begin{aligned}
	\dd \xi_t &= (1-2\eta_t)\; \dd t\\
	\xi_0(x)&=x\end{aligned}\right.
	\quad\quad\text{and}\quad\quad
	\left\{\begin{aligned}
	\dd \eta_t &= 0\;\dd t \\
	\eta_0(x)&=g(x).
	\end{aligned}
	\right. 
	\]
	Hence we obtain
	\begin{equation}\label{eq:deterministic-CE}
		\eta_t(x)=g(x),\qquad 
		\xi_t(x)=x+\int\limits_0^t 1-2g(x)\,\dd s = x+t-2g(x)t.
	\end{equation}
	Then the solution to (\ref{eq:LWR_H=0}) is given by $u(x,t)=g(\xi_t^{-1}(x))$, where $\xi_t^{-1}$ denotes the inverse function of $\xi_t$.
\end{lemma}

\begin{example}\label{EX1:deterministic_g(x)=x}
	Consider the following partial differential equation (PDE) for the Lighthill-Whitham-Richards model on $[0,1]$
	\begin{equation}\label{traffic_deterministic_x}
		\dd u =-(1-2u)\cdot u_x\,\dd t, \qquad u(x,0)=g(x)
	\end{equation}
	Due to Lemma \ref{LE:H=0_CE} we obtain with the initial condition $g(x)=1-x$: 
	\begin{equation}
		\begin{aligned}
			\xi_t(x)&=x-t+2xt,\\
			\eta_t(x)&= 1-x,
		\end{aligned}\qquad\text{hence }\quad  \xi_t^{-1}(x)=\frac{x+t}{1+2t}.
	\end{equation}
	
	Thus the solution of the above PDE (\ref{traffic_deterministic_x}) is given by
	\begin{equation}\label{sol_EX1}
		u(x,t)=\frac{1-x+t}{1+2t}.
	\end{equation}
	
	If we change the initial condition to be $g(x)= 1-x^2$, we obtain 
	\begin{equation}
		\begin{aligned}
			\xi_t(x)&=x-t+2x^2t,\\
			\eta_t(x)&= 1-x^2,
		\end{aligned}\qquad \text{and hence }\qquad \xi_t^{-1}(x)=\frac{\sqrt{1+8t^2+8tx}-1}{4t}.
	\end{equation}
	Thus the corresponding solution is equal to\begin{equation}\label{sol_EX2}
		u(x,t)=\begin{cases}1-\frac{(\sqrt{1+8t^2+8tx}-1)^2}{16t^2}, &\text{for} \; t\neq 0 \\ 1-x^2, &\text{for}\; t=0.
		\end{cases}
	\end{equation}
	One can easily verify that (\ref{sol_EX2}) solves indeed the PDE (\ref{traffic_deterministic_x}).
\end{example}

The main advantage of this method is the precise expression of a solution to a PDE - provided that the corresponding initial condition $g(x)$ and coefficient functions are explicitly given. Due to this fact and for a better comparison between the deterministic and stochastic case we present a collection of solutions in Appendix \ref{app1}.

Along these paths in space-time the solution is constant. In the case of the traffic problem and under the considered initial conditions the characteristics never cross each other which means that no shocks appear and hence the solutions are global.
As written in the introduction we will study the perturbed case (\ref{traffic-eq+H}) for $H(u)\neq 0$ and for $M_t$ to be the standard Brownian motion as well as the so called geometric Brownian motion defined in the following way:
\begin{definition}\label{def_gBM}
	A stochastic process $S_t$, $t\geq 0$, is said to be a \textbf{geometric Brownian motion} if it satisfies
	\begin{equation}\label{SDE_gBM}
		\dd S_t = \mu S_t\,\dd t + \sigma S_t \,\dd W_t= S_t\,\dd W_t
	\end{equation}
	where $W_t$ is a Brownian motion.
	Hence the geometric Brownian motion without drift is given by \[S_t:=\exp{\left(-\frac{t}{2}+W_t\right)}.\]
\end{definition}
\begin{definition}\label{F-notation}
	Let $W_t$ be a standard one-dimensional Brownian motion on a complete, separable probability space $(\Omega,\mathscr{F},{P},\mathscr{F}_t)$ with right-continuous filtration $(\mathscr{F}_t)_{t>0}$. Then we define for any smooth function $H(x,u,p,t)$, $x,u,p\in[0,1]$, $t\in[0,\bT]$, for $0<\bT<\infty$ the following integral expression
	\begin{equation*}
		\int\limits_0^t F(x,u,p,\circ \dd s):=\begin{cases}
			\quad\int\limits_0^t -(1-2u)p\;\dd s + \int\limits_0^t H(x,u,p,s)\circ \dd W_s,\\
			\quad\int\limits_0^t -(1-2u)p\;\dd s + \int\limits_0^t H(x,u,p,s)\circ \dd S_s, 
		\end{cases}
	\end{equation*}
	The integrals are given in the sense of Stratonovich.
\end{definition}

Based on these definitions we are able to apply a heuristic approach of the so called method of stochastic characteristics. Since we consider partial differential equations with perturbations by (geometric) Brownian motions we get an $\omega$\,-\,dependence in the solutions. The idea of the method is nearly the same as before: now we fix $\omega\in\Omega$ and transform a stochastic partial differential equation (SPDE) into a system of stochastic differential equations (SDEs), solve it and determine the solution to the original SPDE by using stopping times. Hence the precisely determined solutions are given for almost all $\omega$ and all space and time variables $(x,t)$ up to a certain stopping time denoted by $\sigma(x)$. In contrast to the deterministic case we will introduce in the following the method of stochastic characteristics in a more detailed way. Based on Definition \ref{F-notation} a perturbed Lighthill-Whitham-Richards model (\ref{traffic-eq+H}) is equivalent to the Cauchy problem
\begin{equation}\label{(SPDE)}
	\left\{\begin{aligned}\dd u &= F(x,u,u_x,\circ \dd t),\\
		u&=g ~~\text{on}~~ \Gamma:=\{x\in[0,1]\times [0,\bT]\,|\, x=(x_1,t), t=0\}.\end{aligned}\right.
\end{equation}

Therefore the solution to equation (\ref{(SPDE)}) is denoted by $u(x,t,\omega)$, but for short notation we only write $u(x,t)$. Suppose $u$ is a solution to (\ref{(SPDE)}) and at least one-times continuously differentiable with respect to space and time for fixed $\omega\in\Omega$. Furthermore, we assume that there exists a curve $\xi_s(r)$ which maps the point $r\in\Gamma$ to a point of a neighborhood in $\Gamma$ at time $s$. Additionally, we assume $\xi_0(x)=x$ for all $x\in[0,1]$ as the initial condition. Due to these assumptions we consider and define the following functions, now for fixed $\omega$, $r\in[0,1]$ and $s\in [0,\bT]$: 
\begin{equation}\label{xi-eta-chi}
	\begin{aligned}
		&(\xi_s(r,\omega),s)\\
		&\eta_s(r,\omega):= u(\xi_s(r,\omega),s),\\
		&\chi_s(r,\omega):=u_{\xi_s}(\xi_s(r,\omega),s).
	\end{aligned}
\end{equation}
In the next step we combine (\ref{(SPDE)}) with equations (\ref{xi-eta-chi}) and obtain 
\[\frac{\dd }{\dd t}\left[u(\xi_t(r),t) -u(\xi_0(r),0)- \int\limits_0^t F(\xi_s(r),\eta_s(r),\chi_s(r),\circ \dd s)\right]=0.\]
By similar calculations as in \cite[§ 3.2.1, equation (11)]{Evans} we get
\begin{equation}\label{(SCE)}\tag{SCE}
	\left\{\begin{aligned} \dd {\xi}_t &= -F_{{\chi}_t}(\xi_t,\eta_t,\chi_t,\circ\dd t),  \\
		\dd {\eta}_t&=F(\xi_t,\eta_t,\chi_t,\circ\dd t) -{\chi}_t\cdot F_{{\chi}_t}(\xi_t,\eta_t,\chi_t,\circ\dd t)  \\
		\dd {\chi}_t &= F_{{\xi}_t}(\xi_t,\eta_t,\chi_t,\circ\dd t)+F_{{\eta}_t}(\xi_t,\eta_t,\chi_t,\circ\dd t){\chi}_t.
	\end{aligned}\right.
\end{equation}

The above stochastic differential equations (\ref{(SCE)}) are called \textbf{stochastic characteristic equations}. Given a point $x\in [0,1]$ and assuming that there exist unique solutions to (\ref{(SCE)}) starting from $x$ at time $t=0$, these solutions solve the corresponding integral equation with initial function $g$:
\begin{align*}
	\xi_t(x)&=x -\int\limits_0^t F_{{\chi}_s}({\xi}_s(x),{\eta}_s(x),{\chi}_s(x),\circ\dd s)\\
	\eta_t(x)&=g(x) -\int\limits_0^t {\chi}_s\cdot F_{{\chi}_s}({\xi}_s(x),{\eta}_s(x),{\chi}_s(x),\circ\dd s)+ \int\limits_0^t F({\xi}_s(x),{\eta}_s(x),{\chi}_s(x),\circ\dd s)\\
	\chi_t(x)&=g_x(x) +\int\limits_0^t F_{{\xi}_s}({\xi}_s(x),{\eta}_s(x),{\chi}_s(x),\circ\dd s)+\int\limits_0^t  F_{{\eta}_s}({\xi}_s(x),{\eta}_s(x),{\chi}_s(x),\circ\dd s){\chi}_s.
\end{align*}

Let us assume that the solutions $({\xi}_t(x),{\eta}_t(x),{\chi}_t(x))$  exist up to a stopping time ${T}(x)$. As mentioned above we have to work on different stopping times based on the following definition.
\begin{definition}\label{stopping-times}
	Let $T(x)$ be the explosion time of the solutions $({\xi}_t,{\eta}_t,{\chi}_t)$ which means e.g. in the case of $\xi_t(x)$ if $$\lim\limits_{t\nearrow T(x)} |\xi_t(x)|=\infty.$$ Then we define for all $x,y\in [0,1]$ the stopping times
	\begin{align*}
		\tau_{\textnormal{inv}}(x)&:=\inf\{t>0\,|\,\det D {\xi}_t(x)=0\},\\
		\tau(x)&:=\tau_{\textnormal{inv}}(x)\wedge {T}(x), \\
		\sigma(y)&:=\inf\{t>0\,|\, y\notin  {\xi}_t(\{x\in[0,1] \,|\,\tau(x)>t\})\},
	\end{align*}
	where $D {\xi}_t$ denotes the Jacobian matrix.
\end{definition}

Now let the inverse process $\xi_t^{-1}$ of ${\xi}_t$ exist up to some stopping time $\sigma(x)$. Then we define for almost all $\omega$ and for all $(x,t)$ with $ t<\sigma(x,\omega)$ the solution
\begin{equation}\label{u-kandidat}
	u(x,t):={\eta}_t({\xi}_t^{-1}(x)).
\end{equation}

Detailed derivations and introductions can be found in \cite[Chapter 3]{diss_NM}. Now we are able to solve different SPDEs concerning the Lighthill-Whitham-Richards model by using the heuristic approach .

\section{Application \& Representation}
Based on the flow rate function $H(u)$ and the continuity equation the most natural choice of the drift term is $H(u)=u-u^2$. In a first step we perturb the Lighthill-Whitham-Richards model by a standard Brownian motion. Hence we consider
\begin{equation}\label{EX1_BM}
	\left\{\begin{aligned}
		\dd u &=-(1-2u)\cdot u_x\,\dd t - (1-2u)\cdot u_x\circ \dd W_t\\
		u(x,0)&=g(x).
	\end{aligned}\right.
\end{equation}
By using the heuristic approach one can show that the corresponding stochastic characteristic equations are given for almost all $\omega$ and all $(x,t)$ up to a stopping time $\sigma(x)$ by
\begin{equation}
	\left\{\begin{aligned}
		\dd \xi_t&=(1-2\eta_t)\,\dd t + (1-2\eta_t)\circ \dd W_t\\
		\xi_0(x)&=x \end{aligned}\right.\quad\quad \text{and}\quad\quad
	\left\{ \begin{aligned}\dd \eta_t&=0 \,\dd t + 0 \circ \dd W_t\\
		\eta_0(x)&=g(x)
	\end{aligned}\right.
\end{equation}
Due to the linearity in the space derivative $u_x$ the solution $\eta_t(x)=g(x)$ is always valid. Therefore we receive the solution
\[\xi_t(x)=x+(1-2g(x))(t+W_t).\]
At this point we compare the characteristics in the deterministic case with the corresponding perturbed one, see Figure 1. As the initial condition we use here $g(x)=1-x$.

\begin{figure}[h]
	\begin{center}
		\hspace*{-3cm}	\begin{minipage}{.4\linewidth}
			\includegraphics[width=1.5\linewidth]{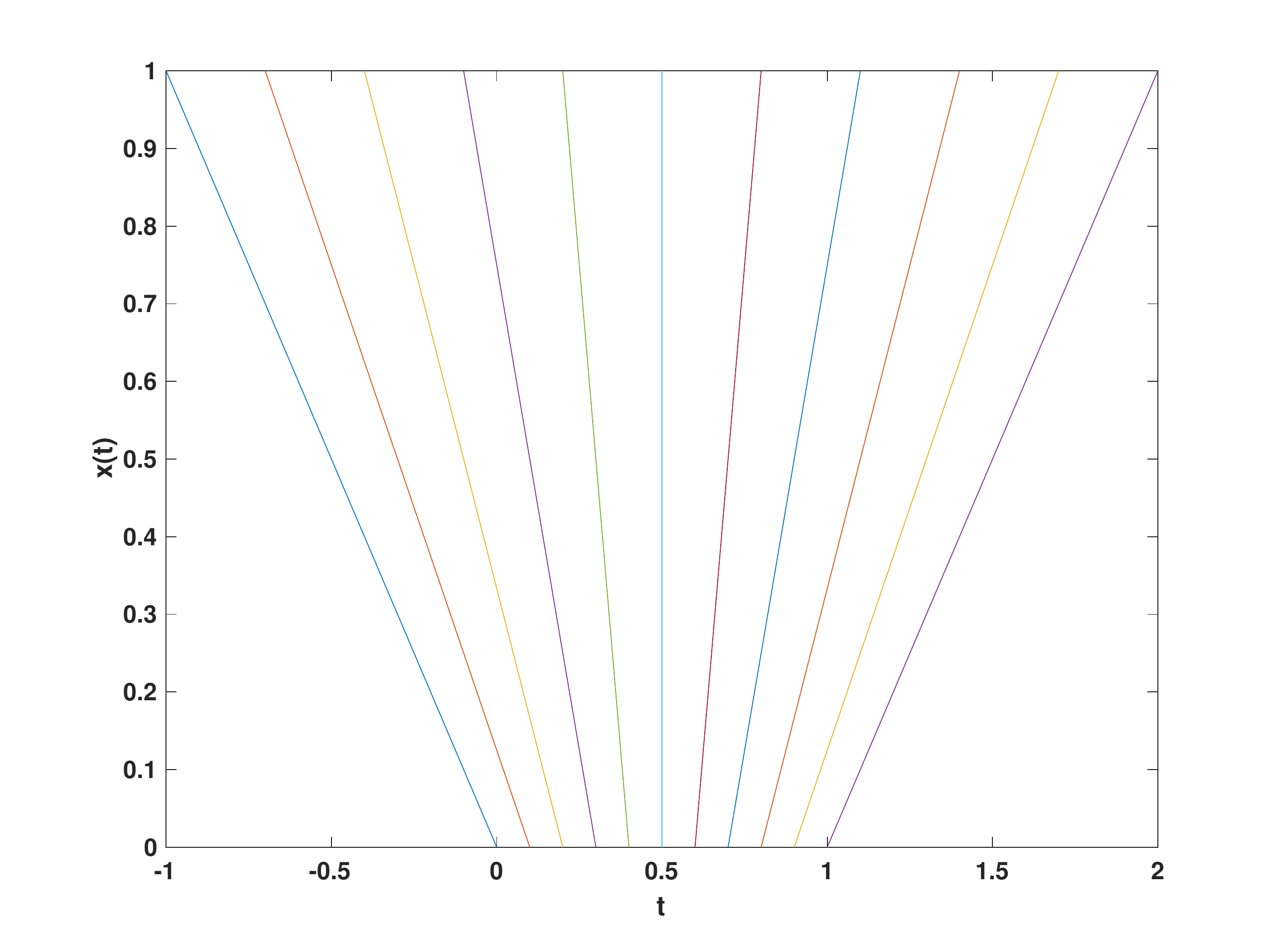} \hfill
		\end{minipage}\hspace{2.5cm}
		\begin{minipage}{.4\linewidth}
			\includegraphics[width=1.5\linewidth]{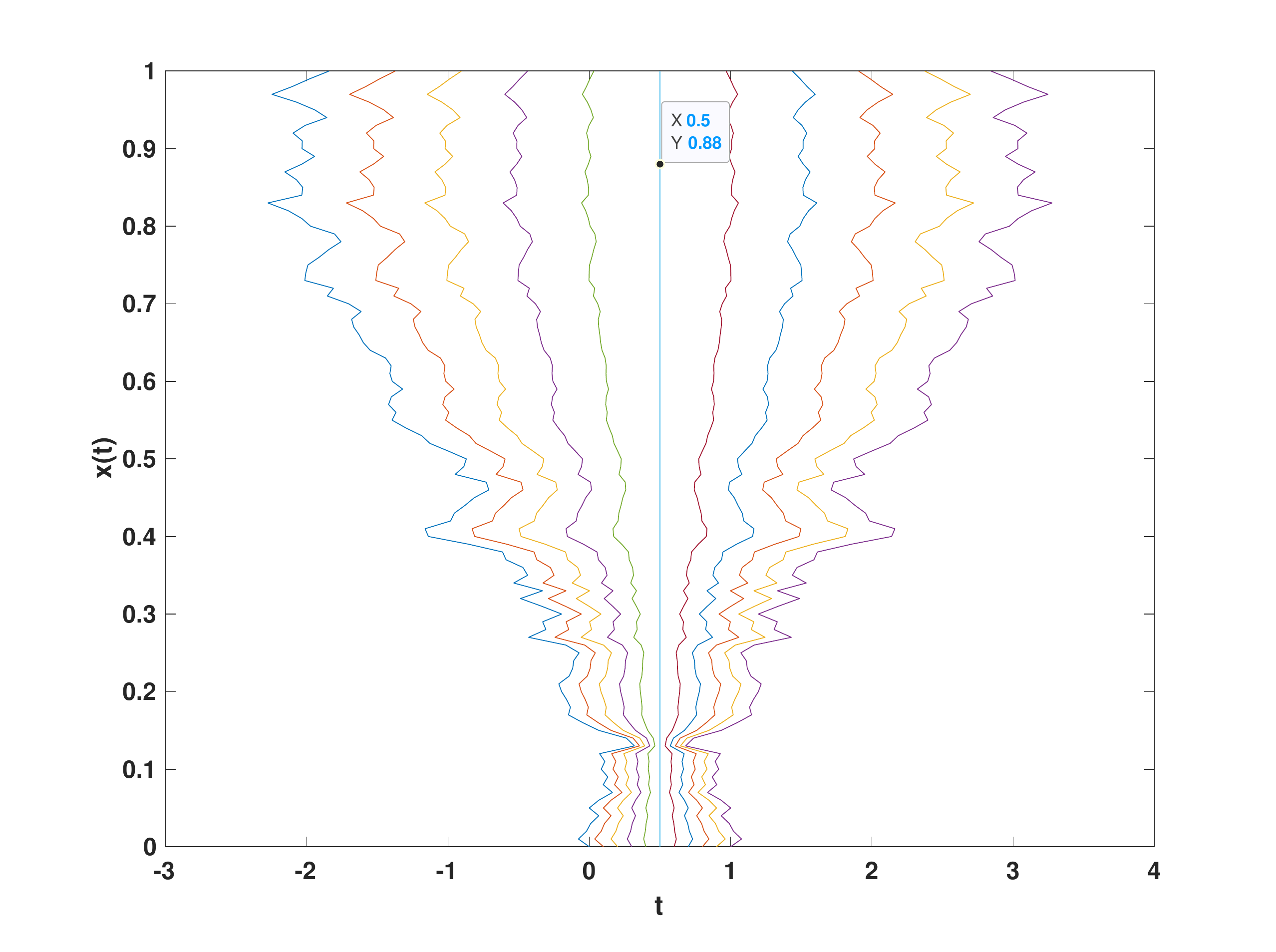}
			\label{fig:EX2_characteristics}
		\end{minipage}
		\caption{Characteristics for the Lighthill-Whitham-Richards model with and without stochastic perturbation}
		\label{fig:EX1_characteristics}
	\end{center}
\end{figure}

In this case of $g(x)=1-x$ there exists obviously a process $\xi_t^{-1}$, such that the inverse property is fulfilled for almost all $\omega$ and all $(x,t)$ up to stopping time $\sigma(x)$, i.e.
$$
\xi_t^{-1}(x)=\frac{x+t+W_t}{1+2t+2W_t}.
$$ 
The solution to the considered SPDE (\ref{EX1_BM}) is given for almost all $\omega$ and all $(x,t)$ up to stopping time $\sigma(x)$ by
\begin{equation}\label{sol_EX1_BM}
	u(x,t)=\frac{1-x+t+W_t}{1+2t+2W_t},
\end{equation}
which looks similar to the deterministic solution (\ref{sol_EX1}). Due to the explicit expression of the solution we are able to visualize a sample path easily, see Figure 2.
\begin{figure}[h]
	\begin{center}
		\includegraphics[width=.6\linewidth]{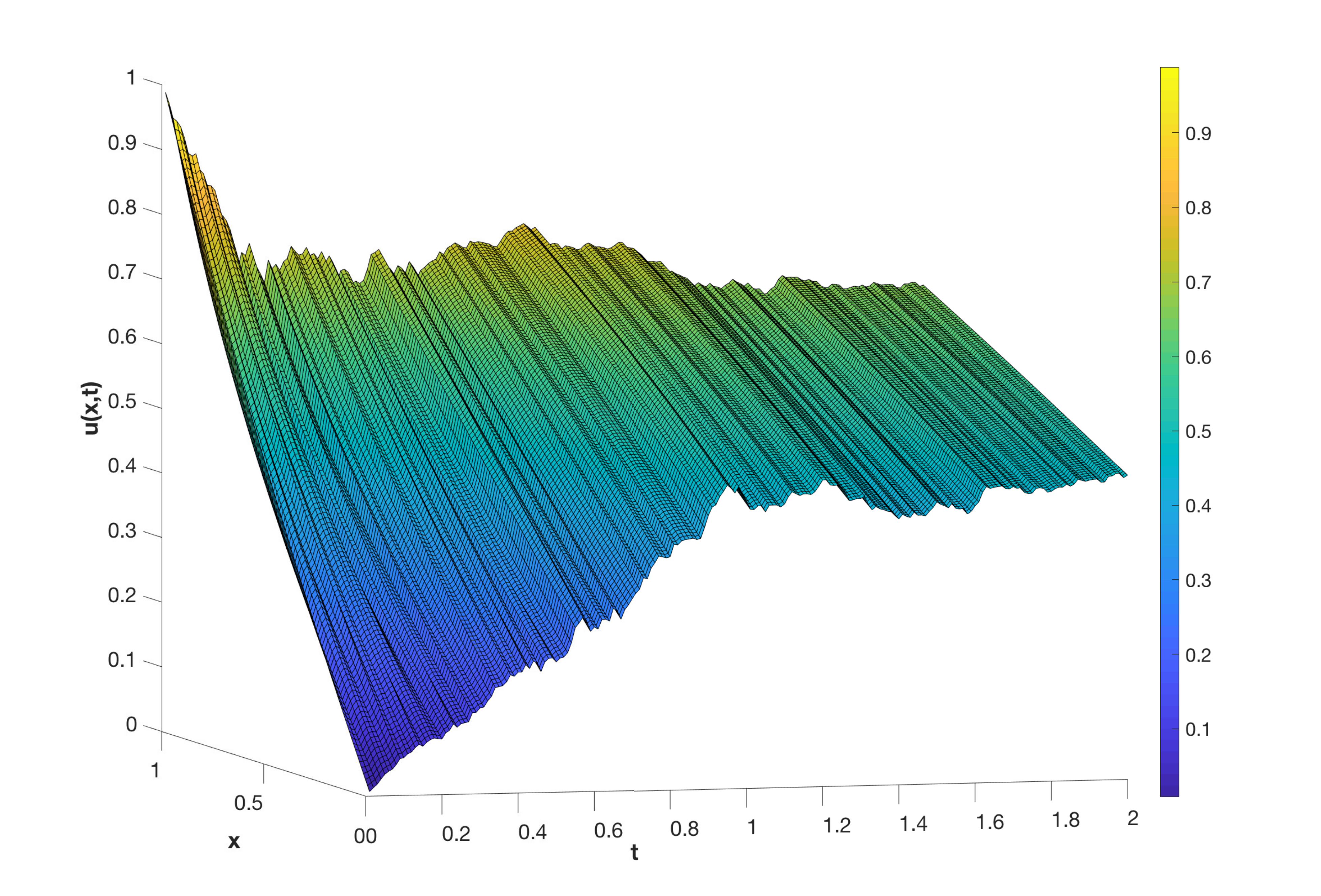}
	\end{center}
	\caption{Sample path of the Lighthill-Whitham-Richards model with initial condition $1-x$ perturbed by the term $- (1-2u)\cdot u_x\circ \dd W_t$. Here we took a sample path without doubling points due to stochasticity.} 
\end{figure}
As introduced in Definition \ref{stopping-times} the stopping time can be determined explicitly in this example by
\begin{align*}
	\sigma(x)&=\inf\Big\{t>0\,\Big|\, \frac{x+t+W_t}{1+2t+2W_t}\notin [0,1]\Big\}\wedge \inf\Big\{t>0\,\Big|\, 1+2t+2W_t=0\Big\} \\
	&=\inf\Big\{t>0\,\Big|\, \frac{x+t+W_t}{1+2t+2W_t}\notin [0,1]\Big\}\wedge \infty.
\end{align*}

The perturbation by a geometric Brownian motion as given in Definition \ref{def_gBM} is in this case straightforward. According to Definition \ref{F-notation} we practically can replace the Brownian motion $W_t$ by $\exp(-t/2+W_t)-1$. Let us consider
\begin{equation}
	\left\{\begin{aligned}
		\dd u &=-(1-2u)\cdot u_x\,\dd t - (1-2u)\cdot u_x\circ \dd [\exp(-t/2+W_t)],\\
		u(x,0)&=1-x^2.
	\end{aligned}\right.
\end{equation}
By an application of the heuristic method of stochastic characteristics we finally get the precise solution for almost all $\omega$ and $(x,t)$ up to a stopping time $\sigma(x)$ by
\begin{equation}\label{sol_EX2_gBM}
	u(x,t)=\begin{cases} 1-\frac{(\sqrt{8(t+e^{(-t/2+W_t)}-1)(t+x+e^{(-t/2+W_t)}-1)+1}-1)^2}{16(e^{(-t/2+W_t)}+t-1)^2}, &\text{if}\; t\neq 0\\
		1-x^2, &\text{if}\; t=0, \end{cases} 
\end{equation}

where we can use the classical l'Hospital argument. The corresponding stopping time is equal to
\begin{align*}
	\sigma(x)&=\inf\Big\{t>0\,\Big|\, \frac{\sqrt{8(e^{(-t/2+W_t)}+t-1)(t+x+e^{(-t/2+W_t)}-1)+1}-1}{4(e^{(-t/2+W_t)}+t-1)}\notin [0,1]\Big\}\\
	&\quad\wedge \inf\Big\{t>0\,\Big|\, \frac{1}{\sqrt{8(t+e^{(-t/2+W_t)}-1)(t+x+e^{(-t/2+W_t)}-1)+1}}=0\Big\}\\
	&=\inf\Big\{t>0\,\Big|\, \frac{\sqrt{8(e^{(-t/2+W_t)}+t-1)(t+x+e^{(-t/2+W_t)}-1)+1}-1}{4(e^{(-t/2+W_t)}+t-1)}\notin [0,1]\Big\} \wedge \infty
\end{align*}
In Figure 3 we display one sample path  with initial condition $1-x^2$ perturbed by the term $- (1-2u)\cdot u_x\circ \dd [\exp(-t/2+W_t)]$.
\begin{figure}[h]
	\begin{center}
		\includegraphics[width=.6\linewidth]{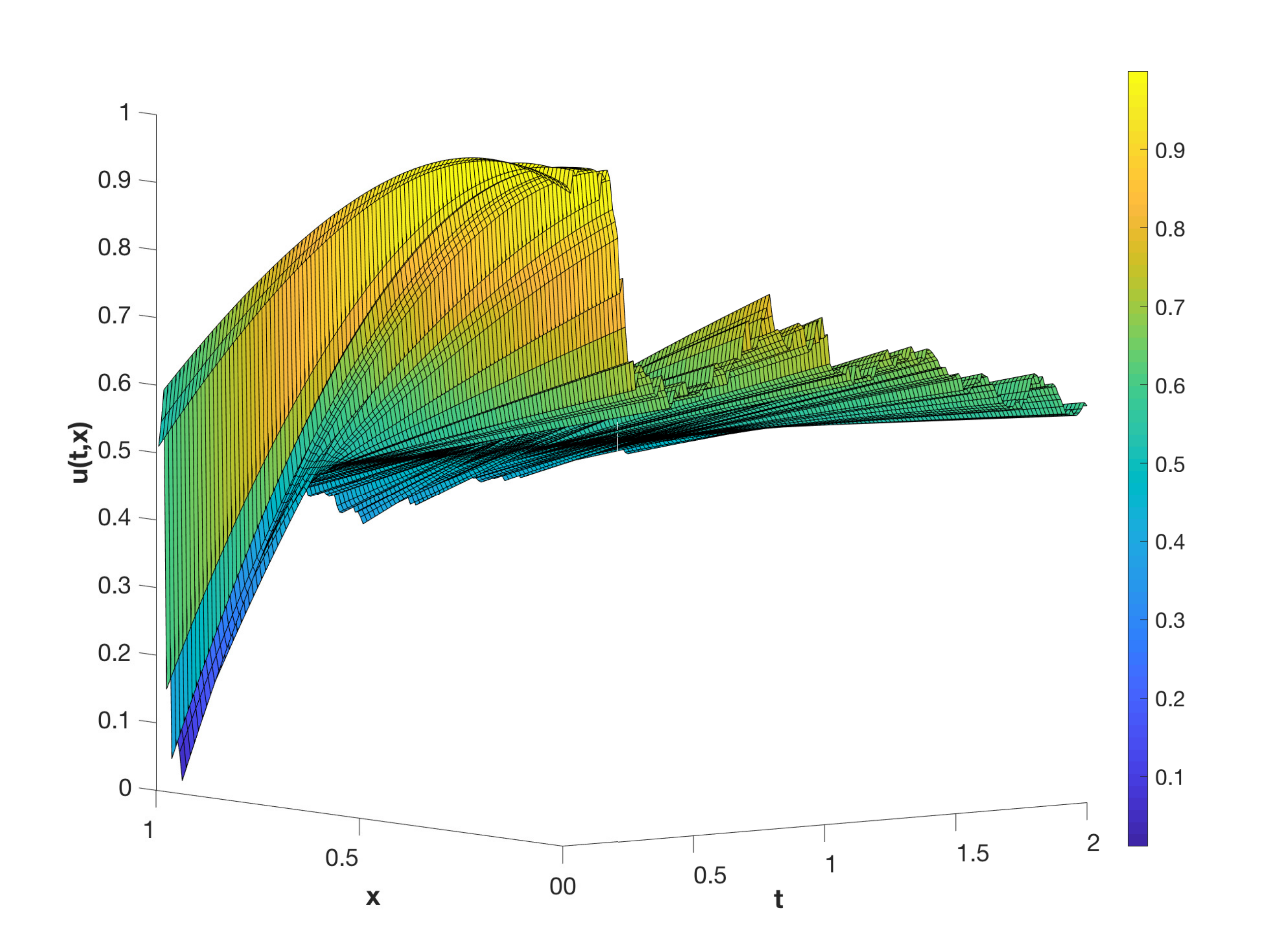}
	\end{center}
	\caption{Sample path of the Lighthill-Whitham-Richards model with initial condition $1-x^2$ perturbed by the term $- (1-2u)\cdot u_x\circ \dd [\exp(-t/2+W_t)]$. Here we took a sample path wwithout doubling points due to stochasticity.} 
\end{figure}

Due to this heuristic approach we have to verify that the equations (\ref{sol_EX1_BM}) as well as (\ref{sol_EX2_gBM}) really solve the underlying problems. For the sake of simplicity these necessary but lengthy calculation can be found in the Appendix \ref{app3}.

For reader's convenience we add some other examples in Appendix \ref{app2} with precise expressions of solutions and different choices of $H(u)$, but which may not rigorously fit the Lighthill-Whitham-Richards model.

\section{Conclusions and Discussion}
The method of stochastic characteristics can be used effectively to solve a stochastic perturbed Lighthill-Whitham-Richards model. The solutions are explicitly given up to a stopping time in closed form. Numerical simulations based on these models can hence been implemented straightforward. However one has to be careful, that the intersection of characteristics due to stochastic perturbation can lead to solutions which are only defined on a smaller time interval than the non-perturbed ones. 
On the other hand, it may be also possible, that the stochastic perturbations increase the time interval where solutions are defined. An example for a solution which is ill-defined due to intersecting characteristics can be seen in Figure 4.

\begin{figure}[h]
	\begin{center}
		\includegraphics[width=.6\linewidth]{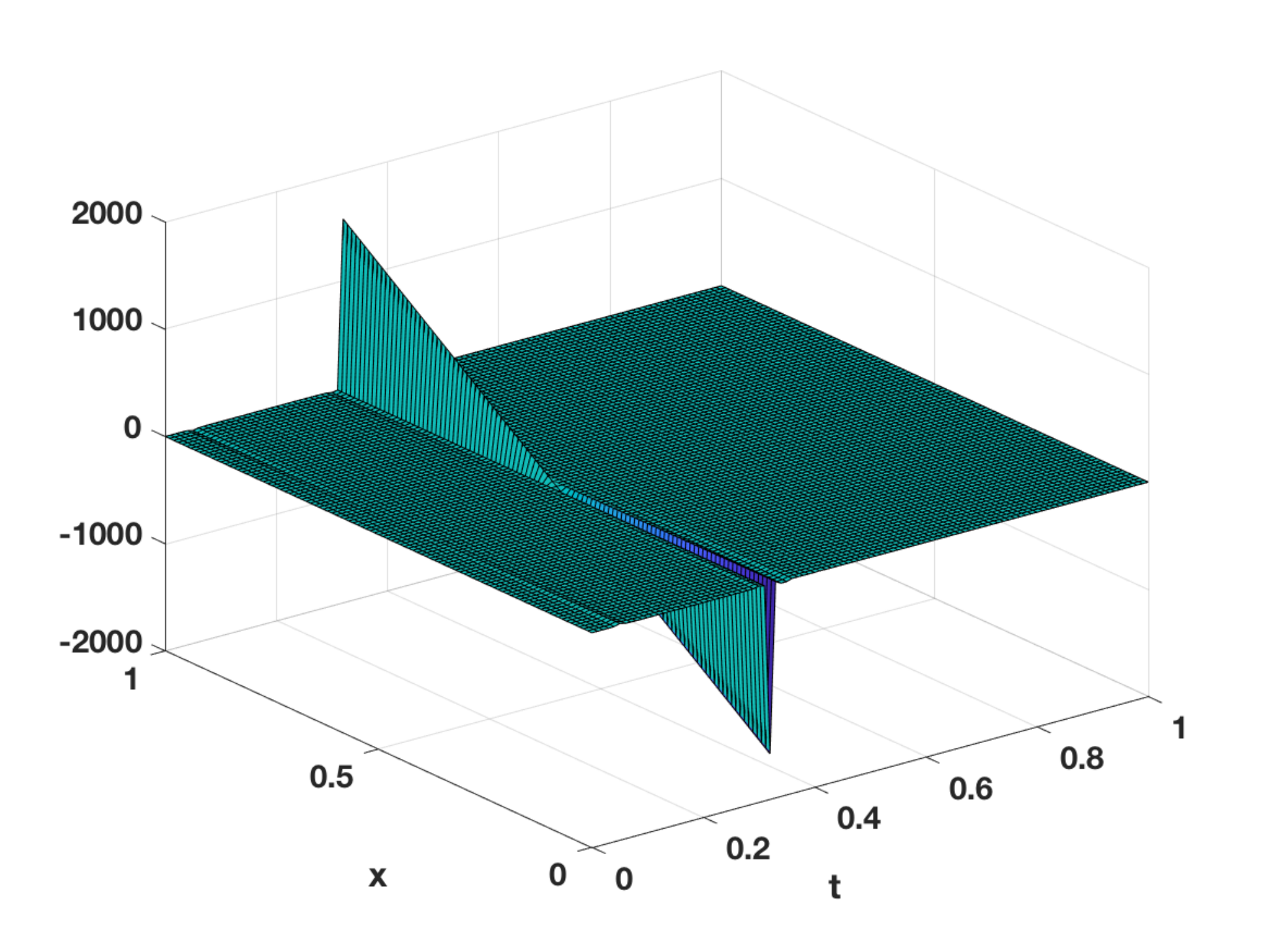}
	\end{center}
	\caption{Sample path of the Lighthill-Whitham-Richards model with initial condition $1-x$ perturbed by the term $- (1-2u)\cdot u_x\circ \dd W_t$. The solution is just defined on a small time interval.} 
\end{figure}

A collection for different examples of stochastic perturbations can be found in Appendix \ref{app2}. Note that with the considered perturbations measurement errors can be modeled effectively. This could be of high interest for more complicated traffic flow models.

\subsection*{Acknowledgement}
We would like to thank for the hospitality of CIMA-UMa at Madeira Math Encounters 2019 in Funchal. Without this meeting point, this publication would not have been possible. Also, we would like to extend our particular thanks to the CRC 1283 \textit{Taming uncertainty and profiting from randomness and low regularity in analysis, stochastics and their applications} for support and generosity in attending conferences and workshops.

\appendix
\section{Collection of examples in the deterministic case}\label{app1}

For reader's convenience the authors itemize the corresponding solutions to the deterministic Lighthill-Whitham-Richards model (\ref{eq:LWR_H=0}) for different initial functions $g(x)$. Based on the model a couple of initial conditions are possible apart from Example \ref{EX1:deterministic_g(x)=x} with $g(x)=1-x$. The opposite to the above case is $g(x)=x$, i.e. the road at position $x=0$ has empty density but at $x=1$ there is e.g. a tailback or a red light, hence the initial density is maximal. We also want to consider a quadratic form by $g(x)=1-x^2$ or $g(x)=x-x^2$, which coincide with the behaviour of the drift part. Analogous calculations yield the following solutions:

\begin{table}[H]
	\begin{tabular}{|l|l}
		\hline
		for & \multicolumn{1}{l|}{solution to
			$\Bigg\{\begin{aligned}
			\dd u &=-(1-2u)\cdot u_x\,\dd t\\
			u(x,0)&= g(x),\quad x\in[0,1]
			\end{aligned}$
		} \\ \hline
		$g(x)=x$   &    $u(x,t)=\frac{x-t}{1-2t}$, $t\neq \frac{1}{2}$                      \\ \cline{1-1}
		$g(x)=1-x^2$    & $u(x,t)=1-\frac{(\sqrt{8t^2+8tx+1}-1)^2}{16t^2}$, $t\neq 0$                         \\ \cline{1-1}
		$g(x)=x-x^2$    & $u(x,t)=\frac{\sqrt{-4t^2+t(8x-4)+1}+2t-1}{4t}-\frac{(\sqrt{-4t^2+t(8x-4)+1}+2t-1)^2}{16t^2}$, $t\neq 0$                         \\ \cline{1-1}
	\end{tabular}
\end{table}

\section{Collection of examples in the stochastic case} \label{app2}

Analogously to the observation in Appendix \ref{app1} we specify the solutions to the perturbed Lighthill-Whitham-Richards model for different choices of the initial function $g(x)$ as well as for different diffusion terms $H(u)$. Taking into account that these might not be model the original traffic flow problem perfectly, the heuristic approach of the method of stochastic characteristics will give explicit solutions. Firstly we perturb the equation by standard Brownian motion.

\begin{itemize}
	\item The solution to the equation
	\begin{equation}
		\left\{\begin{aligned}
			\dd u &=-(1-2u)\cdot u_x\,\dd t + u_x\circ \dd W_t\\
			u(x,0)&=1-x^2
		\end{aligned}\right.
	\end{equation}
	is given for almost all $\omega$ and all $(x,t)$ up to a certain stopping time by
	$$ u(x,t)=1-\frac{(\sqrt{8t(W_t+t+x)+1}-1)^2}{16t^2}
	$$
	\item The solution to the equation
	\begin{equation}
		\left\{\begin{aligned}
			\dd u &=-(1-2u)\cdot u_x\,\dd t + u\circ \dd W_t\\
			u(x,0)&=x
		\end{aligned}\right.
	\end{equation}
	is given for almost all $\omega$ and all $(x,t)$ up to a certain stopping time by
	$$ u(x,t)=\frac{t-x}{2\int\limits_0^t \exp(W_s)\,\dd s -1}
	$$
	\item The solution to the equation
	\begin{equation}
		\left\{\begin{aligned}
			\dd u &=-(1-2u)\cdot u_x\,\dd t + \sqrt{u-u^2}\cdot u_x \circ \dd W_t\\
			u(x,0)&=x
		\end{aligned}\right.
	\end{equation}
	is given for almost all $\omega$ and all $(x,t)$ up to a certain stopping time by
	$$ u(x,t)=\frac{W_t\sqrt{W_t^2+4t^2-4t-4x^2+4x}+W_t^2+4t^2-4xt-2t+2x}{2(1-4t+4t^2+W_t^2)}
	$$
\end{itemize}
Replacing the standard Brownian motion by the geometric Brownian motion without drift we are able to determine also explicit solutions to different SPDEs.

\begin{itemize}
	\item The solution to the equation
	\begin{equation}
		\left\{\begin{aligned}
			\dd u &=-(1-2u)\cdot u_x\,\dd t + u_x\circ \dd [\exp(-t/2+W_t)]\\
			u(x,0)&=x
		\end{aligned}\right.
	\end{equation}
	is given for almost all $\omega$ and all $(x,t)$ up to a certain stopping time by
	$$ u(x,t)=\frac{1-x+t-\exp(-t/2+W_t)}{2t-1}
	$$
	\item The solution to the equation
	\begin{equation}
		\left\{\begin{aligned}
			\dd u &=-(1-2u)\cdot u_x\,\dd t + u\circ \dd [\exp(-t/2+W_t)]\\
			u(x,0)&=x
		\end{aligned}\right.
	\end{equation}
	is given for almost all $\omega$ and all $(x,t)$ up to a certain stopping time by
	$$ u(x,t)=\frac{\e(x-t)}{\e-2\int\limits_0^t \exp(\exp(-s/2+W_s))\,\dd s}
	$$
	\item The solution to the equation
	\begin{equation}
		\left\{\begin{aligned}
			\dd u &=-(1-2u)\cdot u_x\,\dd t + \sqrt{u-u^2}\cdot u_x \circ \dd [\exp(-t/2+W_t)]\\
			u(x,0)&=x
		\end{aligned}\right.
	\end{equation}
	is given for almost all $\omega$ and all $(x,t)$ up to a certain stopping time by
	\begin{align*}
		u(x,t)&=\big(\exp(-t/2+W_t)-1\big)\\
		&\quad\quad \cdot\Big(\sqrt{4t^2-4t-4x^2+4x+\exp(-t+2W_t)-2\exp(-t/2+W_t)+1}\\
		&\quad\quad\quad\quad+4t^2-4xt-2t+2x+\exp(-t+2W_t)-2\exp(-t/2+W_t)+1\Big)\\
		&\quad\quad \cdot\Big(2\big(4t^2-4t+\exp(-t+2W_t)-2\exp(-t/2+W_t)+2\big)\Big)^{-1}
	\end{align*}
\end{itemize}

Formally all given solutions need a verification, similarly to the proofs in Appendix \ref{app3}. But this should not be part of this manuscript.
\section{Calculation and Proofs}\label{app3}
\textbf{Claim:} (\ref{sol_EX1_BM}) solves the stochastic partial differential equation (\ref{EX1_BM})\\
\begin{proof}
	In a first step we determine the partial derivatives $u_t$ and $u_x$ by using $\frac{\circ\dd W_t}{\dd t}=\dot{W_t}$. We obtain
	\begin{align*}
		\frac{\dd u}{\dd t} &=\frac{\dd }{\dd t}\left[\frac{1-x+t+W_t}{1+2t+2W_t}\right]\\
		&=\frac{(1+2t+2W_t)(1+\dot{W_t})-(1-x+t+W_t)(2+2\dot{W_t})}{(1+2t+2W_t)^2}
	\end{align*}
	and
	\begin{align*}
		\frac{\dd u}{\dd x} &=\frac{\dd }{\dd x}\left[\frac{1-x+t+W_t}{1+2t+2W_t}\right]\\
		&=-\frac{1}{(1+2t+2W_t)}
	\end{align*}
	Finally we have to verify that $u_t+(1-2u)u_x+(1-2u)u_x\dot{W_t}=0$.
	\begin{align*}
		u_t&+(1-2u)u_x+(1-2u)u_x\dot{W_t}\\
		&=\frac{(1+2t+2W_t)(1+\dot{W_t})-(1-x+t+W_t)(2+2\dot{W_t})}{(1+2t+2W_t)^2}+\frac{(1-2x)}{(1+2t+2W_t)^2}+\frac{(1-2x)\dot{W_t}}{(1+2t+2W_t)^2}\\
		&=\frac{1+\dot{W_t}+2t+2t\dot{W_t}+2W_t+2W_t\dot{W_t}-(2+2\dot{W_t}-2x-2x\dot{W_t}+2t+2t\dot{W_t}+2W_t+2W_t\dot{W_t})}{(1+2t+2W_t)^2}\\
		&\quad\quad +\frac{1-2x+\dot{W_t}-2x\dot{W_t}}{(1+2t+2W_t)^2}\\
		&=\frac{1+\dot{W_t}+2t+2t\dot{W_t}+2W_t+2W_t\dot{W_t}-2-2\dot{W_t}+2x+2x\dot{W_t}-2t-2t\dot{W_t}-2W_t-2W_t\dot{W_t}}{(1+2t+2W_t)^2}
		\\
		&\quad\quad +\frac{1-2x+\dot{W_t}-2x\dot{W_t}}{(1+2t+2W_t)^2}\\
		&=\frac{\dot{W_t}-2\dot{W_t}+2x-2x+\dot{W_t}}{(1+2t+2W_t)^2}\\
		&=0
	\end{align*}
\end{proof}

\end{document}